# The Smallest Singular Value of a Shifted Random Matrix

Xiaoyu Dong xydong@umich.edu

# THE SMALLEST SINGULAR VALUE OF A SHIFTED RANDOM MATRIX

XIAOYU DONG

ABSTRACT. Let $R_n$ be a $n \times n$ random matrix with i.i.d. subgaussian entries. Let $M$ be a $n \times n$ deterministic matrix with norm $\|M\| \leq n^\gamma$ where $1/2 < \gamma < 1$. The goal of this paper is to give a general estimate of the smallest singular value of the sum $R_n + M$, which improves an earlier result of Tao and Vu.

## 1. INTRODUCTION

### 1.1. Smallest Singular Value of Random Square Matrices.
For an $n \times n$ matrix $A$, recall that the smallest singular value is defined as $s_{\min}(A) = s_n(A) = \inf_{x \in \mathbb{S}^{n-1}} \|Ax\|_2$, where $\mathbb{S}^{n-1}$ is the unit sphere and $\|\cdot\|_2$ is the Euclidean norm.

The early results of the smallest singular value of random square matrices are summarized in the survey [9]. Among the early results, a celebrated work [1] deals with the gaussian random matrices. It shows that, if a matrix $A$ is a $n \times n$ random matrix whose entries are independent standard normal random variables, then the following non-asymptotic bound holds for all $n$:

$$\mathbb{P}(s_{\min}(A) \leq \epsilon n^{-1/2}) \leq \epsilon, \epsilon \geq 0$$

Later Rudelson and Vershynin [8] proved that same estimate holds for general subgaussian random matrices up to a constant coefficient in front of $\epsilon$ and an exponentially small term that accounts for the singularity probability. In detail, they showed that if a matrix $A$ is a $n \times n$ random matrix whose entries are independent and identically distributed subgaussian random variables with mean zero and unit variance, then

$$\mathbb{P}(s_{\min}(A) \leq \epsilon n^{-1/2}) \leq C\epsilon + c^n, \epsilon \geq 0$$

where $c \in (0, 1)$ is an absolute constant. This result is an improvement of the first polynomial bound on $s_{\min}(A)$ for general random matrices that was obtained in [6].

After the work of [8], many new researches appeared and were devoted to relaxing the assumptions in [8] and achieving greater generalities. A non-exhaustive list includes the work by Livshyts, Tikhomirov, and Vershynin [4], the work by Tikhomirov [14], and the work by Farrell and Vershynin [2].

### 1.2. Smoothed Analysis.
Many researches showed that the smoothed analysis of the smallest singular value of random matrices is also a very important topic. Spielman and Teng [12] introduced the smoothed analysis of algorithms, and used it to explain why the simplex algorithm is often successful in practice. In their context, the running time of the algorithm is related to the condition number of some input matrix $M$. They argue that, even if the input matrix $M$ is ill conditioned, the computer will actually work with some small perturbation $M + R$ where $R$ is a random matrix, and $M + R$ is usually well conditioned, making the algorithm successful. In order to make the above argument complete, estimating the smallest singular value of $M + R$ becomes important. With this motivation, Sankar, Spielman, and Teng [11] showed that if $R$ is a random matrix with i.i.d. standard Gaussian entries, then for any deterministic matrix $M \in \mathbb{R}^{n \times n}$, we have

$$\mathbb{P}(s_{\min}(R + M) \leq \epsilon) \leq C\epsilon\sqrt{n}, \epsilon \geq 0$$

Later, Tao and Vu [13] relaxed the assumptions on the distribution of the entries of the random matrix, and obtained the following result.

**Theorem 1.1** (Theorem 3.2. in [13]). *Let $x$ be a random variable with mean zero and bounded second moment, and let $\gamma \geq 1/2, B \geq 0$ be constants. Then there is a constant $c_{\text{Th}\,1.1.1}(x, \gamma, B)$ depending on $x, \gamma, B$ such that the following holds. Let $R_n$ be the random matrix of size $n$ whose entries are iid copies of $x$, let $M$ be a deterministic matrix satisfying $\|M\| \leq n^\gamma$, and let $A_n := M + R_n$. Then*

$$\mathbb{P}(s_n(A_n) \leq n^{-(2B+1)\gamma}) \leq c_{\text{Th}\,1.1.1}(x, \gamma, B)(n^{-B+o(1)} + \mathbb{P}(\|R_n\| \geq n^\gamma))$$

In addition, Tao and Vu [13] also showed that the dependence on $\gamma$ is necessary here, which is different from what happens the Gaussian noise case.

Recently, under the additional assumption on the $K$-rank of the deterministic matrix $M$, Jain and Sawhney [3] derived a new result. They define the $K$-rank of the matrix $M$ as the number $r$ such that $M$ has exactly $r$ singular values of size at least $K\sqrt{n}$. They show that if $R$ is a random matrix of size $n$ whose entries are iid copies of a centered subgaussian random variable with unit variance, then, for any deterministic $M$ with $K$-rank at most $(1 - \eta)n$, we have

$$\mathbb{P}(s_{\min}(R + M) \leq \epsilon n^{-1/2}) \leq C\epsilon\sqrt{n} + c^n, \epsilon \geq 0$$

where the constants $c \in (0, 1)$ and $C > 0$ depend on $K$ and $\eta$.

Nevertheless, if we consider general deterministic matrix $M$ only with control of its norm but without any other extra restrictions, then [13] is still the best result so far.

The main objective of this paper is to try to improve the result in [13].

**1.3. Main Result.** Our result is the following theorem.

**Theorem 1.2** (Invertibility on the Whole Sphere). *Let $R_n(\omega)$ be an $n \times n$ matrix whose entries are independent copies of a centered subgaussian real random variable of unit variance. Let $M$ be an arbitrary deterministic matrix such that $\|M\| \leqslant n^\gamma$. Let $A_n(\omega) = R_n(\omega) + M$. Then there exist constants $c_{\text{Th}\,1.2.1} > 0, c_{\text{Th}\,1.2.2} > 0, c_{\text{Th}\,1.2.3} > 0, n_{\text{Th}\,1.2.1} \in \mathbb{N}$ such that for any $\gamma$ with $1/2 < \gamma < 1$, and $n > n_{\text{Th}\,1.2.1}$, we have $\mathbb{P}(s_n(A_n) \leqslant \epsilon n^{-\gamma}) \leq c_{\text{Th}\,1.2.2}(n^{\gamma-1/2}\epsilon + \exp(-c_{\text{Th}\,1.2.1}n^{2-2\gamma}))$.*

We compare our result with the result of [13]. In our setting, since $R_n$ is a subgaussian random matrix, we have $\mathbb{P}(\|R_n\| \geq n^{1/2}) \leq \exp(-c_{\text{Th}\,1.1.2}n)$ where $c_{\text{Th}\,1.1.2}$ is some constant.

Notice first that as the constant $c_{\text{Th}\,1.1.1}$ depends on $B$, the smallest probability one can obtain in Theorem 1.1 is of a polynomial type, i.e., $O(n^{-B})$ for a fixed $B$. In contrast to it, Theorem 1.2 yields a probability of the exponential type. This is especially important when one considers the probability of singularity of a matrix $A_n$ which corresponds to $\varepsilon = 0$. Now, let us compare the bounds on the smallest singular value for the same probability guarantee. To this end, choose $B$ such that $n^{\gamma-1/2}\epsilon = n^{-B}$. It follows that $n^{-(2B+1)\gamma} = n^{(2\gamma-2)\gamma}\epsilon^{2\gamma}$. Then we have the following variant of [13, Theorem 3.2.].

**Theorem 1.3** (Variant of Theorem 3.2. in [13]). *Let $x$ be a centered subgaussian random variable, and let $\gamma \geq 1/2, \epsilon \geq 0$ be constants. Then there is a constant $c_{\text{Th}\,1.3.1}(\gamma, \frac{\log(\varepsilon)}{\log(n)})$ depending on $\gamma, \frac{\log(\varepsilon)}{\log(n)}$ such that the following holds. Let $R_n$ be the random matrix of size $n$* *whose entries are iid copies of x, let M be a deterministic matrix satisfying $\|M\| \leq n^\gamma$, and let $A_n := M + R_n$. Then*

$$\mathbb{P}(s_n(A_n) \leq n^{(2\gamma-2)\gamma} \epsilon^{2\gamma}) \leq c_{\text{Th }1.3.1}(\gamma, \frac{\log(\varepsilon)}{\log(n)})(n^{\gamma-1/2+o(1)}\epsilon + \exp(-c_{\text{Th }1.3.2}n))$$

Ignoring the fact that the constant $c_{\text{Th }1.3.1}(\gamma, \frac{\log(\varepsilon)}{\log(n)})$ in Theorem 1.3 depends on $\gamma$ and $\epsilon$, we see that the right hand side of Theorem 1.2 and Variant of [13, Theorem 3.2.] are mainly the same up to the term $o(1)$. So, if in the left hand side we have $n^{(2\gamma-2)\gamma}\epsilon^{2\gamma} \leq \epsilon n^{-\gamma}$, then Theorem 1.2 will give a better probability bound. This condition is equivalent to $\epsilon < n^{-\gamma}$.

Theorem 1.1 applies to a broader class of random matrices than Theorem 1.2. However, when applied to matrices with i.i.d. subgaussian entries, Theorem 1.2 provides several improvements over it. First, in the cases where $\varepsilon < n^{-\gamma}$, Theorem 1.2 yields a better smallest singular value bound for the same probability guarantee. Second, setting $\varepsilon = 0$, we obtain that the probability that the matrix $A_n$ is singular does not exceed $\exp(-c_{\text{Th }1.2.1}n^{2-2\gamma})$, which cannot be achieved using Theorem 1.1. Third, the constants in Theorem 1.2 are universal, while those in Theorem 1.1 depend on $B$ and $\gamma$.

1.4. **Overview of the Argument.** We will describe the heuristics of the proof of Theorem 1.2 below. Since $s_{\min}(A) = s_n(A) = \inf_{x \in \mathbb{S}^{n-1}} \|Ax\|_2$, our task is to bound the probability

$$\mathbb{P}(s_n(A_n) \leqslant \epsilon n^{-\gamma})$$

Following [8], we start with the decomposition of the sphere. A vector on $\mathbb{S}^{n-1}$ is called compressible if it can be approximated by a sparse vector, i.e., a vector with relatively small support. For the precise definition of compressible vector, see Definition 4.1. Otherwise, the vector is called incompressible. In this context, we need to establish invertibility for the two classes. So we need to show that with high probability, we have $\inf_{x \in Compressible} \|Ax\|_2 \geq \epsilon n^{-\gamma}$ and $\inf_{x \in Incompressible} \|Ax\|_2 \geq \epsilon n^{-\gamma}$.

The class of compressible vectors has small metric entropy, so we can use the covering argument. For single sparse vector, the estimate in [13, p.2344] shows that, there exist $c_1, c_2 \in (0, 1)$ such that, for any $x \in \mathbb{S}^{n-1}$, $\mathbb{P}(\|A_n x\|_2 < c_1\sqrt{n}) \leq c_2^n$. Using this estimate and the covering argument, we establish that $\inf_{x \in Compressible} \|Ax\|_2 \geq cn^{1/2}$ with high probability where $c$ is a constant.

For the invertibility of incompressible vectors, we use the tool from [8]. Let $X_1(\omega), ..., X_n(\omega)$ denote the column vectors of $A_n(\omega) = R_n(\omega) + M$, and let $H_k$ denote the span of all column vectors except the $k$th. Then, Lemma 3.5 (Invertibility via distance) in [8, p.615] shows that we can bound the probability $\mathbb{P}(\inf_{x \in Incompressible} \|Ax\|_2 \geq \epsilon n^{-\gamma})$ by $\sum_{l=1}^{n} \mathbb{P}(\text{dist}(X_l, H_l) < \epsilon)$. The distance $\text{dist}(X_l, H_l)$ can be estimated by $|\langle Z, X_k \rangle|$ where $Z$ is a random normal vector of $H_k$. For convenience, we consider $k = n$ and then the similar argument will work for all values of $k$. Then, by conditioning on $Z$, the estimation of probability $\mathbb{P}(|\langle Z, X_n \rangle| < \epsilon)$ reduces to the Littlewood-Offord problem, and we therefore use the standard tool, Least Common Denominator (LCD), developed in [10, p.1737] to deal with this problem.

The novelty of our argument mainly lies in the details of the treatment for the incompressible vector. As stated before, we need to estimate the probability $\mathbb{P}(|\langle Z, X_n \rangle| < \epsilon)$ for fixed $Z$, and then take the expectation with respect to $Z$ being a random normal vector of $H_k$. When the LCD of the vector $Z$ is large, the standard argument, e.g. Lemma 3.5, can give desirable estimates of $\mathbb{P}(|\langle Z, X_n \rangle| < \epsilon)$. So, it suffices to bound the probability that $Z$ has small LCD. Let $A'$ be the $(n-1) \times n$ matrix with rows $X_1^T(\omega), ..., X_{n-1}^T(\omega)$. Then the event that $Z$ has small LCD will be represented as the event $\{\omega \in \Omega : \exists x(x \in \mathbb{S}^{n-1}, x \text{ has small LCD}, A'(\omega)x = 0)\}$. Since this event is a statement about all vectors in $\mathbb{S}^{n-1}$ with small LCD, it is natural to use the epsilon net argument, i.e., to construct an epsilon net of all vectors in $\mathbb{S}^{n-1}$ with small LCD, as in [7]. However, in our setting, the matrix $A'$ may have larger norm than what appears in [7] because of the large norm of $M$, so we need a finer net than what is constructed in [7]. In order to overcome this difficulty, we applied a new version of LCD developed in [10, p.1737], the precise definition of which is in Definition 3.1. However, the application of this new LCD introduces some new difficulties, because if we use this new LCD, then the important simple lower bound in [7, p.97; Lemma 6.1.] will not hold. In order to solve this problem, we use the LCD of a truncated part of the vector rather than the whole vector. Now, in order to construct an epsilon net of all vectors in $\mathbb{S}^{n-1}$ with small truncated LCD, we divide each vector into two parts, one with large support but possibly small norm (we call it residual part), and the other with small support but possibly large norm (we call it sparse part). And then, we construct epsilon nets for the collection of the residual parts, and the collection of the sparse parts separately. For the collection of the residual parts, since the truncated part of the vector that we use to evaluate the LCD is exactly the residual part, we can use the definition of LCD to construct the net. More specifically, the LCD of the residual part is small means that if we scale the residual part by a small magnitude, then the resulting vector will be close to an integer point. Also, since the norm of that integer point is close to the norm of the scaled vector, it is bounded by the magnitude of the scaling. Therefore, by collecting the integer points in an Euclidean ball with moderate radius, and rescaling them, we can obtain a net for the collection of the residual parts. For the collection of the sparse parts, we will use a completely different strategy. The sparse part is not related to the truncated LCD of the vector. But the dimension of the sparse part is low, so we can simply use the volumetric argument to construct a net, and the cardinality of this net will still not be harmful for our final estimates because of the low dimension. The combination of those two different strategies will eventually give an epsilon net of all vectors in $\mathbb{S}^{n-1}$ with small truncated LCD.

1.5. **Acknowledgement.** I would like to thank my advisor Prof. Mark Rudelson for his patient guidance, from the construction of the general idea and the argument, to the details in the writing of the draft. Without his help and guidance, this paper would not exist.

## 2. Notation

The following notation and terminology will be used in the paper. We write the unit sphere as $\mathbb{S}^{n-1}$. We denote the operator norm of a matrix $X$ as $\|X\|$. The cardinality of a set $S$ is written as $\text{card}(S)$. The support of a vector $v$ is denoted as $\text{supp}(v)$ which is the set of indices of the nonzero coordinates of $v$. We write the function $\max\{\log(x), 0\}$ as $\log_+(x)$. We denote the open 2-norm ball in $\mathbb{R}^d$ with center $x$ and radius $r$ as $B_2^d(x, r)$, and the corresponding closed ball is called $\bar{B}_2^d(x, r)$. The probability measure is called $\mathbb{P}$, and the symbol $\mathbb{E}$ means taking the expectation with respect to the probability measure. All random matrices in this paper are assumed to be real. According to the standard terminology, a random variable $Z$ is said to be subgaussian if there exists $B > 0$ such that $\mathbb{P}(|Z| > t) \leq 2\exp(-t^2/B^2)$ for all $t > 0$. Throughout the paper, the symbols $c_1, c_2, ...,$ and $Const, Const', ...$ denote the absolute constants. For important constants, we use the symbols such as $c_{\text{Th 1.1.1}}$, $c_{\text{Le 5.1.1}}$, etc., to denote them.

## 3. Small Ball Probability

Least Common Denominator (LCD) and small ball probability are the standard tools to study the smallest singular value of random matrices, so we use a special form of the LCD introduced in [10, p.1737].

**Definition 3.1** (Log Plus Least Common Denominator). Fix $L > 0$. For a vector $v \in \mathbb{R}^n$, the Log Plus Least Common Denominator (LCD) is defined as

$$\mathrm{LCD}_{L,n}(v) = \inf\{\theta > 0 : \mathrm{dist}(\theta \cdot v, \mathbb{Z}^n) < L\sqrt{\max\{\log(\frac{\|\theta \cdot v\|_2}{L}), 0\}}\}$$

In the context where the dimension $n$ is clear, we sometimes use the short notation $\mathrm{LCD}_L$ instead of $\mathrm{LCD}_{L,n}$.

Sometimes, we need to use the LCD not of the original vector, but of a truncated part of the vector, so we make the following definitions, following [15, p.160].

**Definition 3.2** (Truncated Vector). Let $x \in \mathbb{R}^n$. Assume that $k, n$ are integers with $k \leq n$. $I \subset \{1, 2, ..., n\}$ with $\mathrm{card}(I) = k$. Let $\sigma : \{1, 2, ..., k\} \to I$ be the increasing enumeration of $I$. For a vector $v \in \mathbb{R}^n$, we define $\mathrm{TV}_n(x, I) = (x_{\sigma(1)}, x_{\sigma(2)}, ..., x_{\sigma(k)})$.

**Definition 3.3** (Small Part of a Vector). For a vector $v \in \mathbb{R}^n$, assume that we use a permutation $\sigma$ on $\{1, 2, ..., n\}$ to rearrange the coordinates of $v$ in an ordered sequence, i.e., $|v_{\sigma(1)}| \leq |v_{\sigma(2)}| \leq ... \leq |v_{\sigma(n)}|$. Let $I = \{\sigma(1), \sigma(2), ..., \sigma(k)\}$. Then we define $v_{[1:k]} = \mathrm{TV}_n(v, I)$.

We recall the definition of Levy concentration function.

**Definition 3.4** (Levy Concentration). Let $Z$ be a random variable taking values in $\mathbb{R}$. For every $\epsilon > 0$, the Levy concentration of $Z$ is defined as $\mathcal{L}(Z, \epsilon) = \sup_{u \in \mathbb{R}} \mathbb{P}(|Z - u| < \epsilon)$.

**Lemma 3.5** (Small Ball Probability by Truncated LCD). Let $\xi_1, ...$ be a sequence of independent and identically distributed subgaussian random variables. Fix $L > 0$. For any $d \in \mathbb{N}$ and $v \in \mathbb{R}^d$, we define $S_v = \sum_{i=1}^{d} v_i \xi_i$. Assume that $k, n$ are integers with $k \leq n$. Then there exists a constant $c_{\mathrm{Le}\, 3.5.1}$, such that for a vector $v \in \mathbb{R}^n$, if $I \subset \{1, 2, ..., n\}$., then we have $\mathcal{L}(S_v, \epsilon) \leq c_{\mathrm{Le}\, 3.5.1} L(\frac{\epsilon}{\|\mathrm{TV}_n(v,I)\|_2} + \frac{1}{\mathrm{LCD}_{L,k}(\mathrm{TV}_n(v,I))})$.

*Proof.* The proof of this lemma mainly follows from the the proof of Proposition 6.9 in [15, p.160].

By conditioning on random variables $\xi_i$ where $i \notin I$, we have

$$\mathcal{L}(S_v, \epsilon) \leq \mathcal{L}(\sum_{i \in I} v_i \xi_i, \epsilon) = \mathcal{L}(S_{\mathrm{TV}_n(v,I)}, \epsilon)$$

where we use the observation that $S_{\mathrm{TV}_n(v,I)}$ and $\sum_{i \in I} v_i \xi_i$ have the same distribution.

Then, by Remark 6.5 in [15, p.159], we have

$$\begin{aligned}
\mathcal{L}(S_v, \epsilon) &\leq \mathcal{L}(S_{\mathrm{TV}_n(v,I)}, \epsilon) \\
&= \mathcal{L}(\frac{S_{\mathrm{TV}_n(v,I)}}{\|\mathrm{TV}_n(v,I)\|_2}, \frac{\epsilon}{\|\mathrm{TV}_n(v,I)\|_2}) \\
&\leq c_{\mathrm{Le}\, 3.5.1} L(\frac{\epsilon}{\|\mathrm{TV}_n(v,I)\|_2} + \frac{1}{\mathrm{LCD}_{L,k}(\mathrm{TV}_n(v,I))})
\end{aligned}$$

□

## 4. Decomposition of the Sphere

The smallest singular value of $A$ is the minimum of $\|Ax\|_2$ over all vectors $x \in \mathbb{S}^{n-1}$. So, in order to estimate the smallest singular value of $A$, we need to deal with all the vectors $x \in \mathbb{S}^{n-1}$. Following previous research, e.g. [7], we will decompose $\mathbb{S}^{n-1}$ into two families, compressible vectors and incompressible vectors, and deal with each family separately.

So, we make the following definition, according to [7].

**Definition 4.1** (Compressible and incompressible vectors). Let $\delta \in (0,1)$ and $\rho \in (0,1)$. A vector $x \in \mathbb{R}^n$ is called sparse if $\text{card}(\text{supp}(x)) \leq \delta n$. A vector $x \in \mathbb{S}^{n-1}$ is called compressible if $x$ is within Euclidean distance $\rho$ from the set of all sparse vectors. A vector $x \in \mathbb{S}^{n-1}$ is called incompressible if it is not compressible. The sets of sparse, compressible and incompressible vectors will be denoted by $Sparse = Sparse(\delta)$, $Comp = Comp(\delta, \rho)$ and $Incomp = Incomp(\delta, \rho)$ respectively.

Recall that, in the notation section, we stated that the constants $c_1, c_2, ...$ are absolute constants. Actually, in the argument below, one may find that some of those constants depend on the parameter $L$ in the definition of LCD, or the parameter $\delta$ in the definition of compressible vectors. However, throughout the paper, those parameters are chosen one time and then fixed for the entire argument. So, in this view, all the constants $c_1, c_2, ...$ are still absolute constants.

## 5. Invertibility for Compressible Vectors

The class of compressible vectors has small metric entropy, suggesting that we can use the covering argument.

**Lemma 5.1** (Invertibility for compressible vectors). *Let $R_n(\omega)$ be an $n \times n$ matrix whose entries are independent copies of a centered subgaussian real random variable of unit variance. Let $M$ be an arbitrary deterministic matrix such that $\|M\| \leq n^\gamma$ with $\frac{1}{2} < \gamma < 1$. Let $A_n(\omega) = R_n(\omega) + M$. Let $K \geq 1$. Then there exist $c_{\text{Le }5.1.1} > 0, c_{\text{Le }5.1.2} > 0, c_{\text{Le }5.1.3} > 0$ that depend only on $K$, and such that, for all $0 < \delta < c_{\text{Le }5.1.1}$, we have*

$$\mathbb{P}(\inf_{x \in Comp(\delta, \delta n^{1/2-\gamma})} \|A_n x\|_2 \leq c_{\text{Le }5.1.2} n^{1/2}, \text{ and } \|R_n\| \leq K n^{1/2}) \leq e^{-c_{\text{Le }5.1.3} n}$$

*Proof.* The most useful property of compressible vectors is that they can be approximated by the sparse vectors.

If the event $\inf_{x \in Comp(\delta, \delta n^{1/2-\gamma})} \|A_n x\|_2 \leq s n^{1/2}$, and $\|R_n\| \leq K n^{1/2}$ happens, where $s$ is a constant to be determined later, then we can take an $x \in Comp(\delta, \delta n^{1/2-\gamma})$ such that $\|A_n x\|_2 \leq s n^{1/2}$. Since $x \in Comp(\delta, \delta n^{1/2-\gamma})$, there exists $z \in Sparse(\delta)$ such that $\|x - z\| \leq \delta n^{1/2-\gamma}$. Since $\|R_n\| \leq K n^{1/2}$, we have $\|A_n\| \leq \|R_n\| + \|M\| \leq K n^\gamma$. Therefore, we have $\|A_n z\|_2 \leq \|A_n x\|_2 + \|A_n(z-x)\|_2 \leq s n^{1/2} + K n^\gamma \cdot \delta n^{1/2-\gamma} \leq (s + K\delta) n^{1/2}$.

We obtain the following relationship

$$\{\omega : \inf_{x \in Comp(c, c n^{1/2-\gamma})} \|A_n(\omega) x\|_2 \leq s n^{1/2}, \text{ and } \|R_n(\omega)\| \leq K n^{1/2}\}$$
$$\subset \{\omega : \inf_{x \in Sparse(\delta)} \|A_n(\omega) x\|_2 \leq (s + K\delta) n^{1/2}, \text{ and } \|R_n(\omega)\| \leq K n^{1/2}\}$$

So we only need to work with sparse vectors, and for sparse vectors we can use the covering argument.

For convenience, we will treat sparse vectors with different support separately. Let $E_I = \{x \in \mathbb{R}^n : \text{supp}(x) \subset I\}$. We have

$$\{\omega : \inf_{x \in Sparse(\delta)} \|A_n(\omega)x\|_2 \leqslant (s + K\delta)n^{1/2}, \text{ and } \|R_n(\omega)\| \leqslant Kn^{1/2}\}$$
$$\subset \bigcup_{I \subset [n], card(I) = \lceil \delta n \rceil} \{\omega : \inf_{x \in E_I} \|A_n(\omega)x\|_2 \leqslant (s + K\delta)n^{1/2}, \text{ and } \|R_n(\omega)\| \leqslant Kn^{1/2}\}$$

Therefore, it suffices to estimate the probability of each event

$$\{\omega : \inf_{x \in E_I} \|A_n(\omega)x\|_2 \leqslant (s + K\delta)n^{1/2}, \text{ and } \|R_n(\omega)\| \leqslant Kn^{1/2}\}$$

for fixed $I$.

For single sparse vector $x \in Sparse(\delta)$, we can use the estimate in [13, p.2344]. By this estimate, there exist $c_1, c_2 \in (0, 1)$ such that, for any $x \in \mathbb{S}^{n-1}$, $\mathbb{P}(\|A_n x\|_2 < c_1 \sqrt{n}) \leq c_2^n$.

Also, it follows from [7, p.89] that there exists a constant $c_3 > 0$ such that $\mathbb{P}(\|R_n(\omega)\| > Kn^{1/2}) < \exp(-c_3 n)$.

Hence, we conclude that there exists a constant $c_4 > 0$ such that, for any single vector $y \in \mathbb{S}^{n-1}$, we have

$$\mathbb{P}(\{\omega : \|A_n(\omega)y\|_2 \leqslant c_1 \sqrt{n}, \text{ and } \|R_n(\omega)\| \leqslant Kn^{1/2}\}) < \exp(-c_4 n)$$

Following [3], we define the set $G_K = \{v \in \mathbb{R}^n : \|v\|_2 \leq 1, \|Mv\|_2 \leq 2K\sqrt{n}\}$, and we want to use the covering argument to estimate the probability of the event

$$\{\omega : \inf_{x \in E_I} \|A_n(\omega)x\|_2 \leqslant (s + K\delta)n^{1/2}, \text{ and } \|R_n(\omega)\| \leqslant Kn^{1/2}\}$$

Let $\epsilon > 0$ to be chosen later.

By standard volumetric argument, for example, in [5], we observe that $G_K \cap E_I$ can be covered by $(c_5/\epsilon)^{\delta n}$ translates of $\epsilon(G_K \cap E_I) = (\epsilon G_K) \cap E_I$.

We choose one point from each translate of $\epsilon(G_K \cap E_I) = (\epsilon G_C) \cap E_I$ in the covering to form a net $\mathcal{N}_I$ of cardinality $(c_5/\epsilon)^{\delta n}$.

If the event $\{\omega : \inf_{x \in E_I} \|A_n(\omega)x\|_2 \leqslant (s + K\delta)n^{1/2}, \text{ and } \|R_n(\omega)\| \leqslant Kn^{1/2}\}$ happens, then there exists $x \in E_I$ such that $\|A_n(\omega)x\|_2 \leqslant (s + K\delta)n^{1/2}$, and $\|R_n(\omega)\|_2 \leqslant Kn^{1/2}$. In this case, if we require $(s + K\delta) < K$, then we must have $x \in G_K \cap E_I$, because otherwise

$$\|A_n(\omega)x\|_2 \geq \|Mv\|_2 - \|R_n(\omega)\|_2 \geq 2K\sqrt{n} - K\sqrt{n} \geq K\sqrt{n} \geq (s + K\delta)\sqrt{n}$$

Therefore, there exists a $y \in \mathcal{N}_I$ such that $(x - y) \in 2((\epsilon G_K) \cap E_I)$. Therefore, we have $\|R_n(x-y)\|_2 \leq \|R_n\| \|(x-y)\|_2 \leq 2K\sqrt{n}\epsilon$ and $\|M(x-y)\|_2 \leq 4\epsilon K\sqrt{n}$, so $\|A_n(x-y)\|_2 \leq 8\epsilon K\sqrt{n}$. If we choose $\epsilon$ such that $8\epsilon K = (s + K\delta)$, the relationship

$$\|A_n y\|_2 \leq \|A_n x\|_2 + \|A(x-y)\|_2 \leq 2(s + K\delta)\sqrt{n}$$

follows from the previous inequalities.

If we choose $s$ and $\delta$ such that $2(s+K\delta) < c_1$, then we obtain

$$\mathbb{P}(\{\omega : \inf_{x \in E_I} \|A_n(\omega)x\|_2 \leqslant (s+K\delta)n^{1/2}, \text{ and } \|R_n(\omega)\| \leqslant Kn^{1/2}\})$$
$$\leq \mathbb{P}(\{\omega : \exists y(y \in \mathcal{N}_I, \|A_n(\omega)y\|_2 \leqslant 2(s+K\delta)n^{1/2}, \text{ and } \|R_n(\omega)\| \leqslant Kn^{1/2})\})$$
$$\leq \sum_{y \in \mathcal{N}_I} \mathbb{P}(\{\omega : \|A_n(\omega)y\|_2 \leqslant 2(s+K\delta)n^{1/2}, \text{ and } \|R_n(\omega)\| \leqslant Kn^{1/2}\})$$
$$\leq \operatorname{card}(\mathcal{N}_I) \exp(-c_4 n)$$
$$\leq (c_5/\epsilon)^{\delta n} \exp(-c_4 n)$$
$$= (8Kc_5/(s+K\delta))^{\delta n} \exp(-c_4 n)$$

Since $(8Kc_5/(s+K\delta))^{\delta} < (8Kc_5/s)^{\delta}$, for sufficiently small $\delta$, we have $(8Kc_5/(s+K\delta))^{\delta} < 1$. As a consequence, we conclude that

$$\mathbb{P}(\{\omega : \inf_{x \in E_I} \|A_n(\omega)x\|_2 \leqslant (s+K\delta)n^{1/2}, \text{ and } \|R_n(\omega)\| \leqslant Kn^{1/2}\})$$
$$\leq \exp(-c_6 n)$$

where $c_6 > 0$ is some constant.

Taking the union bound, we further conclude that

$$\mathbb{P}(\{\omega : \inf_{x \in Comp(c,cn^{1/2-\gamma})} \|A_n(\omega)x\|_2 \leqslant sn^{1/2}, \text{ and } \|R_n(\omega)\| \leqslant Kn^{1/2}\})$$
$$\leq \sum_{I \subset [n], card(I) = \lceil cn \rceil} \mathbb{P}(\{\omega : \inf_{x \in E_I} \|A_n(\omega)x\|_2 \leqslant (s+K\delta)n^{1/2}, \text{ and } \|R_n(\omega)\| \leqslant Kn^{1/2}\})$$
$$\leq \binom{n}{\lceil \delta n \rceil} \exp(-c_6 n)$$
$$\leq \exp(4e\delta \log(e/\delta)n - c_6 n)$$

Therefore, if $\delta$ is small enough, then $\exp(4e\delta \log(e/\delta)n - c_6 n) \leq \exp(-(c_6/2)n)$, and in this case, we have

$$\mathbb{P}(\{\omega : \inf_{x \in Comp(\delta,\delta n^{1/2-\gamma})} \|A_n(\omega)x\|_2 \leqslant sn^{1/2}, \text{ and } \|R_n(\omega)\| \leqslant Kn^{1/2}\})$$
$$\leq \exp(-(c_6/2)n)$$

Now, in summary, the requirement for $s$ and $\delta$ are $(s+K\delta) < K$, $2(s+K\delta) < c_1$, $(8Kc_5/(s+K\delta))^{\delta} < (8Kc_5/s)^{\delta} < 1$, and $4e\delta \log(e/\delta) < c_6/2$.

So, we choose $c_{\text{Le }5.1.2} = s = \min\{K/2, c_1/4\}$, then there exists $c_{\text{Le }5.1.1} > 0$ such that, for all $0 < \delta < c_{\text{Le }5.1.1}$ all the requirements for $s$ and $\delta$ will be satisfied. In this case, we have $\mathbb{P}(\inf_{x \in Comp(\delta, \delta n^{1/2-\gamma})} \|A_n x\|_2 \leqslant c_{\text{Le }5.1.2} n^{1/2}, \text{ and } \|R_n\| \leqslant Kn^{1/2}) \leq e^{-c_{\text{Le }5.1.3}n}$ where $c_{\text{Le }5.1.3} = c_6/2$.

□

## 6. Invertibility for Incompressible Vectors

To study the invertibility for incompressible vectors, we need to use the tool in [8, Lemma 3.5 (Invertibility via distance)]. Let $X_1(\omega), ..., X_n(\omega)$ denote the column vectors of $A_n(\omega) = R_n(\omega) + M$, and let $H_k$ denote the span of all column vectors except the $k$th. Then, briefly speaking [8, Lemma 3.5 (Invertibility via distance)] allows us to bound the probability $\mathbb{P}(\inf_{x \in Incomp(\delta, \delta n^{1/2-\gamma})} \|A_n x\|_2 \leqslant \epsilon \delta n^{-\gamma})$ by $\mathbb{P}(\operatorname{dist}(X_n, H_n) < \epsilon)$. See the proof of Lemma 6.10 for the details of this argument.

Now, the distance $\text{dist}(X_n, H_n)$ can be estimated by $|\langle Z, X_n \rangle|$ where $Z$ is a random normal vector of $H_k$. Then, the probability $\mathbb{P}(|\langle Z, X_n \rangle| < \epsilon)$ can be estimated by the small ball probability where we will use the truncated LCD of the random normal $Z$.

In order to achieve those objectives, we first get a simple lower bound of LCD for incompressible vectors.

**Lemma 6.1** (Simple Lower Bound of LCD for Incompressible Vectors). *Fix $L > 0$. Let $0 < c < 1$. Let $n \in \mathbb{N}$, and let $\frac{1}{2} < \gamma < 1$. Let $v \in Incomp(\delta, \delta n^{1/2-\gamma})$. Then there exist a universal constant $c_{\text{Le }6.1.1} > 0$ and $n_{\text{Le }6.1.1} \in \mathbb{N}$ depending only on $\delta$ such that for all $n \in \mathbb{N}$, $n > n_{\text{Le }6.1.1}$, there exists a $k \in \{1, 2, ..., n\}$ with*

$$n - \frac{\delta n}{\log(n)} \leq k, \quad \text{and} \quad \text{LCD}_L(\frac{v_{[1:k]}}{\|v_{[1:k]}\|_2}) \geq c_{\text{Le }6.1.1} \frac{\sqrt{n}}{\log n}.$$

*Proof.* Let $\sigma$ be a permutation rearranging the absolute values of coordinates of $v$ in the increasing order, i.e., $|v_{\sigma(1)}| \leq |v_{\sigma(2)}| \leq ... \leq |v_{\sigma(n)}|$. Let $\eta > 0$ be a constant to be chosen later. Define

$$\kappa = \lfloor \eta(n/(\log(n))^2) \rfloor,$$

and assume that $n > n_0(\eta)$, where $n_0(\eta)$ is sufficiently large to ensure that $\kappa > 1$. Let $j_1 \in \mathbb{N}$ be the largest number satisfying

$$\kappa j_1 < \frac{cn}{\log n}.$$

Assume for a moment that for all $j \in \{1, ..., j_1\}$, we have

(6.1) $$|v_{\sigma(n-j\kappa)}| < \frac{1}{2}|v_{\sigma(n-(j-1)\kappa)}|.$$

This condition means that the coordinates of the rearrangement $|v_{\sigma(1)}| \leq |v_{\sigma(2)}| \leq ... \leq |v_{\sigma(n)}|$ increase rapidly, which would imply that there are only a few relatively large ones. We will show that this contradicts the assumption that $v \in Incomp(\delta, \delta n^{1/2-\gamma})$. Indeed, iterating inequality (6.1) for $j = 1, \ldots, j_1$, we obtain

$$v_{\sigma(n-j_1\kappa)} < (\frac{1}{2})^{j_1} \leq (\frac{1}{2})^{\delta \log n/\eta},$$

where the last inequality follows from the definitions of $j_1$ and $\kappa$. Choose $\eta$ so that the right hand side of this inequality is less than $\delta n^{-\gamma}$ for all sufficiently large $n$. Let

$$I = \{\sigma(n - j_1\kappa), \sigma(n - j_1\kappa + 1), ..., \sigma(n)\}$$

be the set of the largest coordinates of $v$. Then $|v_i| \leq cn^{-\gamma}$ for any $i \in \{1, 2, ..., n\}\setminus I$, and so

$$\sum_{i \in \{1,2,...,n\}\setminus I} v_i^2 \leq \text{card}(\{1, 2, ..., n\}\setminus I) \cdot (cn^{-\gamma})^2 \leq n(cn^{-\gamma})^2.$$

As $|I| < \frac{\delta n}{\log n}$, this contradicts the assumption that $v \in Incomp(\delta, \delta n^{1/2-\gamma})$.

We proved that (6.1) does not hold for some $j \leq j_1$ which means that there exists $k \in \mathbb{N}$ such that

$$0 \leq n - k < \frac{\delta n}{\log n} \quad \text{and} \quad |v_{\sigma(n-k)}| \geq \frac{1}{2}|v_{\sigma(n-k+\kappa)}|.$$

Let $J = \{\sigma(1), \sigma(2), ..., \sigma(k)\}$ and let

$$w = (v_1 \cdot \mathbf{1}_J(1), \ldots, v_n \cdot \mathbf{1}_J(n))$$

be the truncated vector where the coordinates $j \notin J$ are replaced by zeroes. Then $\|w\|_\infty = |v_{\sigma(n-k+\kappa)}|$ while

$$\|w\|_2^2 \geq \sum_{s=0}^{\kappa} (v_{\sigma(n-k+s)})^2 \geq \frac{1}{4}\kappa \|w\|_\infty^2.$$

According to Proposition 7.4 (Simple lower bound for LCD). in [10, p.1738], we have

$$\mathrm{LCD}_L\left(\frac{w}{\|w\|_2}\right) \geq \frac{1}{2} \cdot \frac{\|w\|_2}{\|w\|_\infty} \geq \frac{\sqrt{\kappa}}{4} \geq c_{\mathrm{Le}\,6.1.1} \frac{\sqrt{n}}{\log n},$$

which completes the proof. □

Using small ball probability, we can obtain the following result on invertibility on a single vector, which we will need to use later.

**Lemma 6.2** (Invertibility on a single vector by LCD of small part). *Let $R_n(\omega)$ be an $n \times n$ matrix whose entries are independent copies of a centered subgaussian real random variable of unit variance. Let $M$ be an arbitrary deterministic matrix such that $\|M\| \leq n^\gamma$ with $\frac{1}{2} < \gamma < 1$. Let $X_1(\omega), ..., X_n(\omega)$ denote the column vectors of $A_n(\omega) = R_n(\omega) + M$, and let $A'_n(\omega)$ be the $(n-1) \times n$ matrix with rows $(X_1(\omega))^T, ..., (X_{n-1}(\omega))^T$. Let $x \in \mathbb{S}^{n-1}$. Then, there exist constants $c_{\mathrm{Le}\,6.2.1}, c_{\mathrm{Le}\,6.2.2}$ such that for any $L > c_{\mathrm{Le}\,6.2.1}$ and*

$$t > \frac{1}{\mathrm{LCD}_L(\frac{x_{[1:k]}}{\|x_{[1:k]}\|_2})}$$

*we have*

$$\mathbb{P}(\|A_n'x\|_2 \leq t \|x_{[1:k]}\|_2 n^{1/2}) < (c_{\mathrm{Le}\,6.2.2} t)^{n-1}$$

*Proof.* The result follows directly from [10, p.1742] and [8, p.628]. □

By Lemma 3.5, by conditioning on $Z$, we know that we have a good estimate of $\mathbb{P}(|\langle Z, X_n\rangle| < \epsilon)$ if the truncated LCD of $Z$ is large. So, we want to get an upper bound for the probability that the truncated LCD of $Z$ is small. In order to do so, we define the level sets of LCD, and estimate the probability that the truncated LCD of $Z$ lies in each level.

**Definition 6.3** (Level set of LCD on Sphere). Assume that $L > 0$, $D > 0$, and $n \in \mathbb{N}$. Let $S_{D,n,L} = \{x \in \mathbb{S}^{n-1} : D \leq \mathrm{LCD}_{L,n}(x) \leq 2D\}$.

When we estimate the probability that the truncated LCD of $Z$ lies in each level, we will mainly use the epsilon net argument, so we first construct an epsilon net on level set of LCD on sphere.

**Lemma 6.4** (Epsilon nets on level set of LCD on Sphere). *Fix $L > 0$. Then for any $D > \frac{L}{2}$, there exists constants $c_{\mathrm{Le}\,6.4.2}$ and $c_{\mathrm{Le}\,6.4.3}$ (not depending on $n$) such that, for any $n \in \mathbb{N}$, there exists a $\frac{8L}{D}\sqrt{\log_+ \frac{2D}{L}}$ net for $S_{D,n,L}$ with the property that the net lies in $S_{D,n,L}$ and the cardinality of the net is at most $(c_{\mathrm{Le}\,6.4.2}(D+\sqrt{n})/\sqrt{n})^n$.*

*Proof.* Consider arbitrary $n \in \mathbb{N}$ and $D > 0$. Now, we will construct a net on $S_{D,n,L}$ and calculate its cardinality. For an arbitrary point $x \in S_{D,n,L}$, we will show that $x$ could be approximated by $\frac{p}{\|p\|_2}$ where $p$ is a integer point. For an arbitrary point $x \in S_{D,n,L}$, we know that $\mathrm{LCD}_L(x) \leq 2D$. By Definition 3.1 and continuity of $L\sqrt{\log_+ \frac{\|\lambda \cdot x\|}{L}}$ in $\lambda$, there exists $p \in \mathbb{Z}^n$ such that

$$\|(\mathrm{LCD}_L(x)) \cdot x - p\|_2 \leq L\sqrt{\log_+ \frac{\|(\mathrm{LCD}_L(x)) \cdot x\|}{L}} = L\sqrt{\log_+ \frac{(\mathrm{LCD}_L(x))}{L}}$$

Now, we have
$$\left\| x - \frac{p}{(\mathrm{LCD}_L(x))} \right\|_2 < \frac{L}{(\mathrm{LCD}_L(x))} \sqrt{\log_+ \frac{(\mathrm{LCD}_L(x))}{L}} \leq 2\frac{L}{D}\sqrt{\log_+ \frac{2D}{L}}$$

Also,
$$\left\| \frac{p}{(\mathrm{LCD}_L(x))} \right\|_2 \leq \|x\|_2 + \left\| x - \frac{p}{(\mathrm{LCD}_L(x))} \right\|_2 < 1 + 2\frac{L}{D}\sqrt{\log_+ \frac{2D}{L}}$$

Therefore, we obtain
$$\left\| x - \frac{p}{\|p\|_2} \right\|_2 = \left\| x - \frac{\frac{p}{(\mathrm{LCD}_L(x))}}{\left\| \frac{p}{(\mathrm{LCD}_L(x))} \right\|_2} \right\|_2$$
$$\leq \left\| x - \frac{p}{(\mathrm{LCD}_L(x))} \right\|_2 + \left\| \frac{p}{(\mathrm{LCD}_L(x))} - \frac{\frac{p}{(\mathrm{LCD}_L(x))}}{\left\| \frac{p}{(\mathrm{LCD}_L(x))} \right\|_2} \right\|_2$$
$$\leq 2\frac{L}{D}\sqrt{\log_+ \frac{2D}{L}} + \left(1 - \frac{1}{\left\| \frac{p}{(\mathrm{LCD}_L(x))} \right\|_2}\right) \left\| \frac{p}{(\mathrm{LCD}_L(x))} \right\|_2$$
$$= 2\frac{L}{D}\sqrt{\log_+ \frac{2D}{L}} + \left\| \frac{p}{(\mathrm{LCD}_L(x))} \right\|_2 - 1$$
$$\leq 2\frac{L}{D}\sqrt{\log_+ \frac{2D}{L}} + (1 + 2\frac{L}{D}\sqrt{\log_+ \frac{2D}{L}}) - 1 = 4\frac{L}{D}\sqrt{\log_+ \frac{2D}{L}}$$

Now, we have shown that $x$ can be approximated by $\frac{p}{\|p\|_2}$ where $p$ is a integer point.

Since there are infinitely many integer points, in order to get a finite net, we need to find some restrictions of the integer point $p$. Actually, by triangle inequality and the continuity results, we have
$$\|p\|_2 \leq \|(\mathrm{LCD}_L(x)) \cdot x - p\|_2 + \|(\mathrm{LCD}_L(x)) \cdot x\|_2$$
$$\leq L\sqrt{\log_+ \frac{(\mathrm{LCD}_L(x))}{L}} + \mathrm{LCD}_L(x)$$
$$\leq L\sqrt{\log_+ \frac{2D}{L}} + 2D$$

We observe that $\frac{L}{D}\sqrt{\log_+ \frac{2D}{L}} \leq c_1$ for some absolute constant $c_1$, since $\lim\limits_{D \to \infty} \frac{L}{D}\sqrt{\log_+ \frac{2D}{L}} = 0$. It follows that $L\sqrt{\log_+ \frac{2D}{L}} \leq c_1 D$, so $\|p\|_2 < (3c_1 + 2)D$.

Summarizing, we have shown that the set $\mathcal{N} = \{\frac{p}{\|p\|_2} : p \in \mathbb{Z}^n \cap B_{\mathbb{R}^n, \|\cdot\|_2}(0, (3c_1+2))\}$ is a $\frac{4L}{D}\sqrt{\log_+ \frac{2D}{L}}$ net for $S_{D,n,L}$.

According to standard volumetric argument, e.g. [5, p.11] and [7, p.103], we know that there exists a constant $c_{\mathrm{Le}\,6.4.2}$ not depending on $n$ and a $\frac{8L}{D}\sqrt{\log_+ \frac{2D}{L}}$ net for $S_{D,n,L}$ such that the net lies in $S_{D,n,L}$ and the cardinality of the net is $(c_{\mathrm{Le}\,6.4.2}(D+\sqrt{n})/\sqrt{n})^n$. □

Since, in our setting, we need to use not just the LCD of the original vector $Z$, but of a truncated part of $Z$, we need to further define level set of LCD of small part.

**Definition 6.5** (Level set of LCD of small part). Assume that $L > 0$, $D > 0$, $n \in \mathbb{N}$, and $k \leq n$. Let $\mathrm{TS}_{D,k,n,L} = \{x \in \mathbb{S}^{n-1} : D \leq \mathrm{LCD}_L(\frac{x_{[1:k]}}{\|x_{[1:k]}\|_2}) \leq 2D\}$.

When we switch between the space of truncated vectors and the space of original unit vectors in $\mathbb{R}^n$, it will be convenient to define the following embedding.

**Definition 6.6.** Let $I \subset \{1, 2, ..., n\}$ with $\mathrm{card}(I) = k$. We define $\iota_I$ to be the embedding from $\mathbb{R}^k$ to $\mathbb{R}^n$, such that we just fill coordinates not in $I$ with zeros. Formally, let $\sigma : \{1, 2, ..., k\} \to I$ be the increasing enumeration of $I$. Then $\iota_I$ is the embedding such that for $y \in \mathbb{R}^k$, we have $\iota_I(y)(\sigma(i)) = y(i)$ and $\iota_I(y)(j) = 0$ if $j \notin I$.

**Example 6.7.** Let $n = 4$, $I = \{2, 4\}$, and $y = (0.4, 0.6)$. Then $\iota_I(y) = (0, 0.4, 0, 0.6)$.

Now, we can estimate the probability that the truncated LCD of $Z$ lies in each level.

**Lemma 6.8** (Probability of Random Normal in Level Set). *Let $R_n(\omega)$ be an $n \times n$ matrix whose entries are independent copies of a centered subgaussian real random variable of unit variance. Let $M$ be an arbitrary deterministic matrix such that $\|M\| \leq n^\gamma$. Fix $L > 1$. Let $0 < \delta < 1$. Let $C > 0$. Let $X_1(\omega), ..., X_n(\omega)$ denote the column vectors of $A_n(\omega) = R_n(\omega) + M$, and let $H_k$ denote the span of all column vectors except the kth, for any $1 \leq k \leq n$. Assume that $Z(\omega)$ is a random unit vector orthogonal to $\mathrm{span}(\{X_1(\omega), ..., X_{n-1}(\omega)\})$ for every $\omega \in \Omega$. Then there exists constants $c_{\mathrm{Le}\,6.8.1} > 0, c_{\mathrm{Le}\,6.8.2} > 0$, $n_{\mathrm{Le}\,6.8.1} \in \mathbb{N}$, such that for any $\gamma$ with $1/2 < \gamma < 1$, $D$ with $\frac{1}{2}C\sqrt{n}/(\log(n)) < D < \exp(c_{\mathrm{Le}\,6.8.1}n^{2-2\gamma})$, $n > n_{\mathrm{Le}\,6.8.1}$, and $k$ with $n - \delta n/\log(n) \leq k \leq n$, we have*

$$\mathbb{P}(Z \in \mathrm{TS}_{D,k,n,L} \cap Incomp(\delta, \delta n^{1/2-\gamma})) < \exp(-c_{\mathrm{Le}\,6.8.2}n)$$

*Proof.* We prove the lemma in the following steps. Let $A'$ be the $(n-1) \times n$ matrix with rows $X_1^T(\omega), ..., X_{n-1}^T(\omega)$.

In order to determine whether the event

$$\{\omega : Z(\omega) \in \mathrm{TS}_{D,k,n,L} \cap Incomp(\delta, \delta n^{1/2-\gamma})\}$$
$$= \{\omega \in \Omega : \exists x (x \in Incomp(\delta, \delta n^{1/2-\gamma}) \cap \mathrm{TS}_{D,k,n,L}, A'(\omega)x = 0)\}$$

happens, we need to check all the vectors in $Incomp(\delta, \delta n^{1/2-\gamma}) \cap \mathrm{TS}_{D,k,n,L}$. Since we have a bound of the probability that $A'(\omega)x$ is small for a single vector, it is natural to consider the epsilon-net argument.

In order to enable the epsilon-net argument to work, we need to first control the norm of the random matrix $A'_n(\omega)$. We write $A'_n(\omega) = R'_n(\omega) + M'$. It follows from Proposition 4.4 in [7, p.89] that there exists absolute constants $c_1 > 1$ and $c_2 > 0$ such that $\mathbb{P}(\|R'_n\| > c_1\sqrt{n}) \leq e^{-c_2 n}$. Therefore, we know that $\mathbb{P}(\|A'_n\| > c_1 n^\gamma) \leq e^{-c_2 n}$.

Let

$$\mathcal{E} = \{\omega \in \Omega : \exists x (x \in Incomp(\delta, \delta n^{1/2-\gamma}) \cap \mathrm{TS}_{D,k,n,L}, A'(\omega)x = 0, \|A'(\omega)\| \leq c_1 n^\gamma)\}$$

It then follows that $\mathbb{P}(Z \in \mathrm{TS}_{D,k,n,L} \cap Incomp(\delta, \delta n^{1/2-\gamma})) \leq P(\mathcal{E}) + \mathbb{P}(\|A'_n\| > c_1 n^\gamma) \leq P(\mathcal{E}) + e^{-c_2 n}$. We already have the tool Lemma 6.4 to construct epsilon net on $\mathrm{S}_{D,k,L}$. In order to use Lemma 6.4, we need to fix the index of the small part of the random normal vector, so we consider the following relationship

$$\mathrm{TS}_{D,k,n,L} \subset \bigcup_{I \subset \{1,2,...,n\}, \mathrm{card}(I)=k} \{x \in \mathbb{S}^{n-1} : D \leq \mathrm{LCD}_{L,k}(\frac{1}{\|\mathrm{TV}_n(x,I)\|_2}\mathrm{TV}_n(x,I)) \leq 2D\}$$

$$= \bigcup_{I \subset \{1,2,...,n\}, \mathrm{card}(I)=k} \{x \in \mathbb{S}^{n-1} : \frac{1}{\|\mathrm{TV}_n(x,I)\|_2}\mathrm{TV}_n(x,I) \in \mathrm{S}_{D,k,L}\}$$

For $I \subset \{1, 2, ..., n\}, \text{card}(I) = k$, let

$$\begin{aligned}\mathcal{E}_I =& \{\omega \in \Omega : \exists x (x \in Incomp(\delta, \delta n^{1/2-\gamma}),\\
& x \in \{x \in \mathbb{S}^{n-1}, D \leqslant \text{LCD}_{L,k}(\frac{1}{\|\text{TV}_n(x,I)\|_2} \text{TV}_n(x, I)) \leqslant 2D\},\\
& A'(\omega)x = 0 \text{ and } \|A'(\omega)\| \leqslant c_1 n^\gamma)\}\\
=& \{\omega \in \Omega : \exists x (x \in Incomp(\delta, \delta n^{1/2-\gamma}),\\
& \frac{1}{\|\text{TV}_n(x,I)\|_2} \text{TV}_n(x, I) \in S_{D,k,L}, A'(\omega)x = 0, \|A'(\omega)\| \leqslant c_1 n^\gamma)\}\end{aligned}$$

Therefore, we have $\mathcal{E} \subset \bigcup_{I \subset \{1,2,...,n\}, \text{card}(I) = k} \mathcal{E}_I$, so $\mathbb{P}(\mathcal{E}) \leq \sum_{I \subset \{1,2,...,n\}, \text{card}(I) = k} \mathbb{P}(\mathcal{E}_I)$.

Now we construct epsilon nets and use points in epsilon net to approximate any vector in $\{x \in \mathbb{S}^{n-1} : \frac{1}{\|\text{TV}_n(x,I)\|_2} \text{TV}_n(x, I) \in S_{D,k,L}\}$.

Let $P_I$ be the projection on coordinates in $I$. With this definition, we observe that $P_I(x) = \iota_I(\text{TV}_n(x, I))$.

Since the description of the set $\{x \in \mathbb{S}^{n-1} : \frac{1}{\|\text{TV}_n(x,I)\|_2} \text{TV}_n(x, I) \in S_{D,k,L}\}$ mainly focuses on the properties of $\text{TV}_n(x, I)$, we decompose the vector $x$ into two parts using the identity

$$x = P_I(x) + P_{\{1,2,..,n\}\setminus I}(x) = \iota_I(\text{TV}_n(x, I)) + \iota_{\{1,2,..,n\}\setminus I}(\text{TV}_n(x, \{1, 2, .., n\}\setminus I))$$

and approximate each part separately.

For the first part $\iota_I(\text{TV}_n(x, I))$, we write

$$\iota_I(\text{TV}_n(x, I)) = \iota_I(\|\text{TV}_n(x, I)\|_2 (\frac{1}{\|\text{TV}_n(x, I)\|_2} \text{TV}_n(x, I)))$$

Since $\frac{1}{\|\text{TV}_n(x,I)\|_2} \text{TV}_n(x, I) \in S_{D,k,L}$, we can use a point in an epsilon net of $S_{D,k,L}$ to approximate it. Therefore, using Lemma 6.4, let $\mathcal{N}_1$ be a $\beta_1 = \frac{8L}{D}\sqrt{\log_+ \frac{2D}{L}}$-net for $S_{D,k,L}$ with the property that the net lies in $S_{D,k,L}$ and the cardinality of the net is $(c_3(D + \sqrt{k})/\sqrt{k})^k$. So, there exists $y \in \mathcal{N}_1$ such that $\|y - \frac{1}{\|\text{TV}_n(x,I)\|_2} \text{TV}_n(x, I)\|_2 < \beta_1$. Also, since $\beta_1 = \frac{8L}{D}\sqrt{\log_+ \frac{2D}{L}}$, and $\frac{1}{2}C\sqrt{n}/(\log(n)) < D$, we know that for large enough $n$, we always have $\beta_1 < 1$.

Since $\|\text{TV}_n(x, I)\|_2 \in [0, 1]$, we can use a point in an epsilon net of $[0, 1]$ to approximate it. Therefore, let $\mathcal{N}_3$ be a $\beta_3 = \frac{1}{4}\delta n^{1/2-\gamma}\beta_1$-net for $[0, 1]$ with the property that the net lies in $[0, 1]$ and the cardinality of the net is $(3/\beta_3)$. So, there exists $\xi \in \mathcal{N}_3$ such that $|\xi - \|\text{TV}_n(x, I)\|_2| < \beta_3$.

For the approximation of second part $\iota_{\{1,2,..,n\}\setminus I}(\text{TV}_n(x, \{1, 2, .., n\}\setminus I))$, we use a finer net in $\bar{B}_2^{n-k}(0, 1)$. The reason is that the dimension of $\bar{B}_2^{n-k}(0, 1)$ is small, and therefore even when the net is finer, the cardinality of the net will still be acceptable. Let $\mathcal{N}_2$ be a $\beta_2 = \frac{1}{4}\delta n^{1/2-\gamma}\beta_1$-net for $\bar{B}_2^{n-k}(0, 1)$ with the property that the net lies in $\bar{B}_2^{n-k}(0, 1)$ and the cardinality of the net is $(3/\beta_2)^{n-k}$ by the standard volumetric estimate. So, there exists $z \in \mathcal{N}_2$ such that $\|z - \text{TV}_n(x, \{1, 2, ..., n\}\setminus I)\|_2 < \beta_2$.

Now we have the final approximation of $x$.

$$\begin{aligned}x =& \iota_I(\|\mathrm{TV}_n(x,I)\|_2 \cdot (\frac{1}{\|\mathrm{TV}_n(x,I)\|_2}\mathrm{TV}_n(x,I) - y)) + \iota_I((\|\mathrm{TV}_n(x,I)\|_2 - \xi) \cdot y)\\ &+ \iota_I(\xi \cdot y) + \iota_{\{1,2,..,n\}\setminus I}(\mathrm{TV}_n(x,\{1,2,..,n\}\setminus I) - z) + \iota_{\{1,2,..,n\}\setminus I}(z)\\ =& \iota_I(\xi \cdot y) + \iota_{\{1,2,..,n\}\setminus I}(z) + \iota_I(\|\mathrm{TV}_n(x,I)\|_2 \cdot (\frac{1}{\|\mathrm{TV}_n(x,I)\|_2}\mathrm{TV}_n(x,I) - y))\\ &+ \iota_I((\|\mathrm{TV}_n(x,I)\|_2 - \xi) \cdot y) + \iota_{\{1,2,..,n\}\setminus I}(\mathrm{TV}_n(x,\{1,2,..,n\}\setminus I) - z)\end{aligned}$$

In the above identity, the term $\iota_I(\xi \cdot y) + \iota_{\{1,2,..,n\}\setminus I}(z)$ is the approximation from the points in the epsilon nets, and the term $\iota_I(\|\mathrm{TV}_n(x,I)\|_2 \cdot (\frac{1}{\|\mathrm{TV}_n(x,I)\|_2}\mathrm{TV}_n(x,I) - y)) + \iota_I((\|\mathrm{TV}_n(x,I)\|_2 - \xi) \cdot y) + \iota_{\{1,2,..,n\}\setminus I}(\mathrm{TV}_n(x,\{1,2,..,n\}\setminus I) - z)$ is the error.

Since $\iota_I(\xi \cdot y) + \iota_{\{1,2,...,n\}\setminus I}(z)$ and $x$ are approximations of each other, we can also approximate $\|A'(\omega)(\iota_I(\xi \cdot y) + \iota_{\{1,2,..,n\}\setminus I}(z))\|_2$ by $\|A'(\omega)x\|_2$, using triangle inequality. If $\|A'(\omega)x\|_2 = 0$, we have

$$\begin{aligned}&\|A'(\omega)(\iota_I(\xi \cdot y) + \iota_{\{1,2,..,n\}\setminus I}(z))\|_2\\ \leq& \|A'(\omega)x\|_2 + \|A'(\omega)(\iota_I(\xi \cdot y) + \iota_{\{1,2,..,n\}\setminus I}(z) - x)\|_2\\ \leq& \|A'(\omega)x\|_2 + \left\|A'(\omega)(\iota_I(\|\mathrm{TV}_n(x,I)\|_2 \cdot (\frac{1}{\|\mathrm{TV}_n(x,I)\|_2}\mathrm{TV}_n(x,I) - y))\right\|_2\\ &+ \|A'(\omega)\iota_I((\|\mathrm{TV}_n(x,I)\|_2 - \xi) \cdot y)\|_2 + \|A'(\omega)\iota_{\{1,2,..,n\}\setminus I}(\mathrm{TV}_n(x,\{1,2,..,n\}\setminus I) - z))\|_2\\ \leq& \|\mathrm{TV}_n(x,I)\|_2 \cdot n^\gamma \cdot \beta_1 + n^\gamma \cdot \beta_3 + n^\gamma \cdot \beta_2\end{aligned}$$

where we use first $\|A'(\omega)x\|_2 = 0$, second $\left\|A'(\omega)(\iota_I(\|\mathrm{TV}_n(x,I)\|_2 \cdot (\frac{1}{\|\mathrm{TV}_n(x,I)\|_2}\mathrm{TV}_n(x,I) - y))\right\|_2 \leq \|\mathrm{TV}_n(x,I)\|_2 \cdot n^\gamma \cdot \beta_1$ because $\|A'(\omega)\| \leq n^\gamma$ and $\left\|y - \frac{1}{\|\mathrm{TV}_n(x,I)\|_2}\mathrm{TV}_n(x,I)\right\|_2 < \beta_1$, third $\|A'(\omega)\iota_I((\|\mathrm{TV}_n(x,I)\|_2 - \xi) \cdot y)\|_2 \leq n^\gamma \cdot \beta_3$ because $|\xi - \|\mathrm{TV}_n(x,I)\|_2| < \beta_3$, and finally $\|A'(\omega)\iota_{\{1,2,..,n\}\setminus I}(\mathrm{TV}_n(x,\{1,2,...,n\}\setminus I) - z))\|_2 \leq n^\gamma \cdot \beta_2$ because $\|z - \mathrm{TV}_n(x,\{1,2,...,n\}\setminus I)\|_2 < \beta_2$.

Since $x \in Incomp(\delta, \delta n^{1/2-\gamma})$, we know that $\|\mathrm{TV}_n(x,I)\|_2 \geq \delta n^{1/2-\gamma}$. Then

$$\xi \geq \|\mathrm{TV}_n(x,I)\|_2 - \beta_3 \geq \frac{3}{4}\delta n^{1/2-\gamma} \geq \beta_3$$

So, $\|\mathrm{TV}_n(x,I)\|_2 \leq \xi + \beta_3 \leq 2\xi$.

Also, we have $\beta_3 \leq \frac{1}{4}\delta n^{1/2-\gamma}\beta_1 \leq \xi\beta_1$, and $\beta_2 \leq \frac{1}{4}\delta n^{1/2-\gamma}\beta_1 \leq \xi\beta_1$.

Therefore, we have the following upper bound for the error of the approximation

$$\begin{aligned}&\|A'(\omega)(\iota_I(\xi \cdot y) + \iota_{\{1,2,..,n\}\setminus I}(z))\|_2\\ \leq& \|\mathrm{TV}_n(x,I)\|_2 \cdot n^\gamma \cdot \beta_1 + n^\gamma \cdot \beta_3 + n^\gamma \cdot \beta_2\\ \leq& 2\xi n^\gamma \cdot \beta_1 + \xi n^\gamma \cdot \beta_1 + \xi n^\gamma \cdot \beta_1\\ =& 4\xi n^\gamma \cdot \beta_1\end{aligned}$$

For single $\xi \in \mathcal{N}_3, y \in \mathcal{N}_1, z \in \mathcal{N}_2$, according to Lemma 6.2, we know that

$$\mathbb{P}(\{\omega : \|A'(\omega)(\iota_I(\xi \cdot y) + \iota_{\{1,2,..,n\}\setminus I}(z))\|_2 < \xi t n^{1/2}\}) < (c_{\mathrm{Le\,6.2.2}}t)^{n-1}$$

Let $t = 4n^{\gamma-1/2} \cdot \beta_1$. Now we can bound the probability of $\mathcal{E}_I$ by the union bound.

$\mathbb{P}(\mathcal{E}_I) \leq \mathbb{P}(\{\omega : \exists \xi \exists y \exists z(\xi \in \mathcal{N}_3, y \in \mathcal{N}_1, z \in \mathcal{N}_2, \|A'(\omega)(\iota_I(\xi \cdot y) + \iota_{\{1,2,..,n\}\setminus I}(z))\|_2 < \xi t n^{1/2})\})$
$\leq \sum_{\xi \in \mathcal{N}_3, y \in \mathcal{N}_1, z \in \mathcal{N}_2} \mathbb{P}(\{\omega : \|A'(\omega)(\iota_I(\xi \cdot y) + \iota_{\{1,2,..,n\}\setminus I}(z))\|_2 < \xi t n^{1/2}\})$

$$\leq \operatorname{card}(\mathcal{N}_1) \operatorname{card}(\mathcal{N}_2) \operatorname{card}(\mathcal{N}_3)(c_{\text{Le 6.2.2}} 4n^{\gamma-1/2} \cdot \beta_1)^{n-1}$$
$$\leq ((c_3(D+\sqrt{k})/\sqrt{k})^k)((3/\beta_2)^{n-k})(3/\beta_3)(c_{\text{Le 6.2.2}} 4n^{\gamma-1/2} \cdot \beta_1)^{n-1}.$$

Now, taking the union bound, we eventaully conclude the estimate

$$\mathbb{P}(\mathcal{E}) \leq \binom{n}{k} \cdot \mathbb{P}(\mathcal{E}_{k,I})$$
$$\leq 2^n \cdot ((c_3(D+\sqrt{(k)}))/\sqrt{(k)})^{(k)})((3/\beta_2)^{n-k})(3/\beta_3)(c_{\text{Le 6.2.2}} 4n^{\gamma-1/2} \cdot \beta_1)^{n-1}$$
$$\leq \exp(c_6 n) \cdot (((D+\sqrt{k})/\sqrt{k})(n^{\gamma-1/2} \cdot \beta_1))^k ((1/\beta_2)^{n-k})(1/\beta_3)(n^{\gamma-1/2} \cdot \beta_1)^{n-k-1}$$

with some constant $c_6$ because we use the fact $\binom{n}{k} \leq 2^n$ and $k$ is close to $n$, i.e., $n/2 \leq k \leq n$.

Now, we will simplify the probability bound. From the expression, intuitively, we see that the terms $(((D+\sqrt{k})/\sqrt{k})(n^{\gamma-1/2} \cdot \beta_1))^k$ should be the main terms and $((1/\beta_2)^{n-k})(1/\beta_3)(n^{\gamma-1/2} \cdot \beta_1)^{n-k-1}$ should be the junk terms because $k$ is close to $n$ and therefore $k$ is much larger than $n-k$.

Let $r > 0$ be a constant to be chosen later. From the simplification, we will show that if the constant $r > 0$ is small enough, and $D < \exp(rn^{2-2\gamma})$, then the desired inequality stated in the Lemma will hold.

Let us first try to simplify the main term. In order to simplify the calculation, let $Const$, $Const'$, ... be constants that depend only on $L$, $c$, and $C$, and can change from line to line. We consider two different cases.

In the first case, we assume that $\frac{1}{2}C\sqrt{n}/(\log(n)) < D < \sqrt{n}$. For large enough $n$, we have $\lceil n/2 \rceil \leq n - \delta n/\log(n)$. This estimate yields that $(D+\sqrt{k})/\sqrt{k} \leq Const$. Therefore, for large enough $n$, we obtain that

$$\text{(main term)}$$
$$= (((D+\sqrt{k})/\sqrt{k})(n^{\gamma-1/2} \cdot \beta_1))^k$$
$$\leq \exp(k((\gamma-1/2)\log(n) + \log(\beta_1)))$$

Also, using the bound $\frac{1}{2}C\sqrt{n}/(\log(n)) < D < \sqrt{n}$, we have

$$\log(\beta_1)$$
$$\leq (\log(8L) - \log(D) + \frac{1}{2}\log(\log_+ \frac{2D}{L}))$$
$$\leq Const - \log(D) + \frac{1}{2}\log(2\log(\sqrt{n}))$$
$$\leq Const - \log(\frac{1}{2}C\sqrt{n}/(\log(n))) + \log\log(n)$$
$$= Const - \frac{1}{2}\log(n) + 2\log\log(n)$$

Replacing $\log(\beta_1)$ by the above upper bound in the main term, we conclude that

$$\text{(main term)}$$
$$= (((D+\sqrt{k})/\sqrt{k})(n^{\gamma-1/2} \cdot \beta_1))^k$$
$$\leq \exp(k(Const'' + \log\log(n) + \log(D) + (\gamma-1)\log(n)$$
$$\quad - \log(D) + \log\log(n)))$$
$$\leq \exp(k(Const'' + 2\log\log(n) + (\gamma-1)\log(n)))$$
$$\leq \exp(-sn)$$

where $s > 0$ could be any constant that we choose, because $\gamma - 1 < 0$.

In the second case, we assume that $D > \sqrt{n}$. This assumption yields that

$$(D + \sqrt{k})/\sqrt{k}$$
$$\leq (Const)(D/(\sqrt{n}))$$

As a result, for large enough $n$, we know that

$$\text{(main term)}$$
$$= (((D + \sqrt{k})/\sqrt{k})(n^{\gamma - 1/2} \cdot \beta_1))^k$$
$$\leq \exp(k(Const' + \log(D) + (\gamma - 1)\log(n) + \log(\beta_1)))$$

Also, using the bound $D < \exp(rn^{2-2\gamma})$, we have

$$\log(\beta_1)$$
$$\leq (\log(8L) - \log(D) + \frac{1}{2}\log(\log_+ \frac{2D}{L}))$$
$$\leq Const - \log(D) + \frac{1}{2}\log(\log(\exp(rn^{2-2\gamma}))) + \log(\frac{2}{L}))$$
$$\leq Const - \log(D) + \frac{1}{2}\log(rn^{2-2\gamma})$$
$$= Const - \log(D) + (1 - \gamma)\log(n) + \frac{1}{2}\log(r)$$

Replacing $\log(\beta_1)$ by the above upper bound in the main term, we obtain

$$\text{(main term)}$$
$$= (((D + \sqrt{k})/\sqrt{k})(n^{\gamma - 1/2} \cdot \beta_1))^k$$
$$\leq \exp(k(Const'' + \log(D) + (\gamma - 1)\log(n)$$
$$- \log(D) + (1 - \gamma)\log(n) + \frac{1}{2}\log(r)))$$
$$\leq \exp(k(Const'' + \frac{1}{2}\log(r)))$$

Therefore, if we choose $r > 0$ small enough, the coefficient $(Const'' + \frac{1}{2}\log(r))$ will become negative and its absolute value will be as large as we want. In this way, for any $s > 0$, the estimate holds

$$\text{(main term)}$$
$$\leq \exp(k(Const'' + \frac{1}{2}\log(r)))$$
$$\leq \exp(-sn)$$

for all large enough $n$.

Then we will evaluate the junk term.

We calculate

(junk term)
$$= \exp((n-k)\log(n^{1/2-\gamma}\beta_1)) + \log(n^{1/2-\gamma}\beta_1)) + (n-k-1)\log(n^{\gamma-1/2} \cdot \beta_1))$$
$$\leq \exp(-2\log(\beta_1) + 2(n-k)(\gamma-1/2)\log(n) + (n-k)(Const) + (Const'))$$
$$\leq \exp(2rn^{2-2\gamma} - 2(1-\gamma)\log(n) - \log(r)$$
$$+ 2(n-k)(\gamma-1/2)\log(n) + (n-k)(Const) + (Const'))$$
$$\leq \exp(2rn^{2-2\gamma} - 2(1-\gamma)\log(n) - \log(r)$$
$$+ 2cn/(\log(n))(\gamma-1/2)\log(n) + (n-k)(Const) + (Const'))$$
$$\leq \exp(2rn^{2-2\gamma} - 2(1-\gamma)\log(n) - \log(r)$$
$$+ 2cn(\gamma-1/2) + (n-k)(Const) + (Const'))$$
$$\leq \exp(c_8 n)$$

where $c_8 > 0$ is a constant, because $n^{2-2\gamma} \leq n$. Recall here that $1/2 < \gamma < 1$.

Combining the main term and the junk term, we conclude that

$$\mathbb{P}(\mathcal{E})$$
$$\leq \exp((c_8 - s)n)$$

Therefore, if we choose $s > 0$ large enough, then the right hand side of the above inequality will be an exponential decay. □

Combining the previous preparations, we can now claim that, with high probability, the truncated LCD of $Z$ is large.

**Lemma 6.9** (Truncated LCD of random normal). *Let $R_n(\omega)$ be an $n \times n$ matrix whose entries are independent copies of a centered subgaussian real random variable of unit variance. Let $M$ be an arbitrary deterministic matrix such that $\|M\| \leq n^\gamma$. Fix $L > 1$. Let $0 < c < 1$. Let $X_1(\omega), ..., X_n(\omega)$ denote the column vectors of $A_n(\omega) = R_n(\omega) + M$, and let $H_k$ denote the span of all column vectors except the kth. Assume that $Z(\omega)$ is a random unit vector such that it is orthogonal to $\mathrm{span}(\{X_1(\omega), ..., X_{n-1}(\omega)\})$ for every $\omega \in \Omega$. Then there exists constants $c_{\mathrm{Le}\,6.9.1} > 0, c_{\mathrm{Le}\,6.9.2} > 0, n_{\mathrm{Le}\,6.9.1} \in \mathbb{N}$, and a function $k(n,v)$ with $n - \frac{\delta n}{\log(n)} \leq k(n,v) \leq n$, such that for any $\gamma$ with $1/2 < \gamma < 1$, and $n > n_{6.901}$, we have $\mathbb{P}(\mathrm{LCD}_L(\frac{Z_{[1:k(n,Z)]}}{\|Z_{[1:k(n,Z)]}\|_2}) < \exp(c_{\mathrm{Le}\,6.9.1} n^{2-2\gamma})) < \exp(-c_{\mathrm{Le}\,6.9.2} n)$.*

*Proof.* We divide the target event $\mathrm{LCD}_L(\frac{Z_{[1:k(n,Z)]}}{\|Z_{[1:k(n,Z)]}\|_2}) < \exp(sn^{2-2\gamma})$, where $s$ is a constant to be chosen later, into two cases, according to whether $Z$ is compressible or incompressible.

For the case when $Z$ is incompressible, we use Lemma 6.8. In order to use Lemma 6.8, we need to choose the number $k$ with $n - \delta n/\log(n) \leq k \leq n$ appropriately.

According to Lemma 6.1, there is constant $c_{\mathrm{Le}\,6.1.1} > 0$, such that we can choose a function $k(n,v)$ defined for any $n, v$ satisfying $v \in Incomp(\delta, \delta n^{1/2-\gamma})$ such that

$$n - \frac{\delta n}{\log(n)} \leq k(n,v), \quad \text{and} \quad \mathrm{LCD}_L(\frac{v_{[1:k(n,v)]}}{\|v_{[1:k(n,v)]}\|_2}) \geq c_{\mathrm{Le}\,6.1.1}\frac{\sqrt{n}}{\log n}.$$

for all $v \in Incomp(\delta, \delta n^{1/2-\gamma})$ and large enough $n$.

For $n, v$ satisfying $v \in Comp(\delta, \delta n^{1/2-\gamma})$, we also define $k(n,v) = \lceil n - \frac{cn}{\log(n)} \rceil$.

Now, according to Lemma 6.8, there exists constants $c_{\text{Le } 6.8.1} > 0, c_{\text{Le } 6.8.2} > 0, n_{\text{Le } 6.8.1} \in \mathbb{N}$ such that for any $\gamma$ with $1/2 < \gamma < 1$, $D$ with $\frac{1}{2}c_{\text{Le } 6.1.1}\sqrt{n}/(\log(n)) < D < \exp(c_{\text{Le } 6.8.1}n^{2-2\gamma})$, $n > n_{\text{Le } 6.8.1}$, we have $\mathbb{P}(Z \in \text{TS}_{D,k(n,Z),n,L} \cap Incomp(\delta, \delta n^{1/2-\gamma})) < \exp(-c_{\text{Le } 6.8.2}n)$.

Let $p_D = \mathbb{P}(Z \in \text{TS}_{D,k(n,Z),n,L} \cap Incomp(\delta, \delta n^{1/2-\gamma}))$. It follows directly from the definition of $k(n,v)$ that $p_D = 0$ if $D < c_{\text{Le } 6.1.1}\frac{\sqrt{n}}{\log n}$.

Hence we have

$$\mathbb{P}(\text{LCD}_L(\frac{Z_{[1:k(n,Z)]}}{\|Z_{[1:k(n,Z)]}\|_2}) < \exp(c_{\text{Le } 6.8.1}n^{2-2\gamma}), Z \in Incomp(\delta, \delta n^{1/2-\gamma}))$$
$$\leq \sum_{D=2^i, c_{\text{Le } 6.1.1}\frac{\sqrt{n}}{\log n} < 2^i < \exp(c_{\text{Le } 6.8.1}n^{2-2\gamma})} p_D$$
$$\leq \sum_{\log(c_{\text{Le } 6.1.1}\frac{\sqrt{n}}{\log n})/\log(2) < i < c_{\text{Le } 6.8.1}n^{2-2\gamma}/\log(2)} p_{2^i}$$
$$\leq ((c_{\text{Le } 6.8.1}n^{2-2\gamma})/(\log(2)))\exp(-c_{\text{Le } 6.8.2}n)$$
$$\leq \exp(-c_8 n)$$

where $c_8 > 0$ is some constant.

Now, we will deal with the case when $Z$ is compressible.

We write $A'_n(\omega) = R'_n(\omega) + M'$. It follows from [7, p.89] that there exists constants $c_4 > 0$ and $c_5 > 0$ such that $\mathbb{P}(\|A'_n\| > c_4 n^\gamma) \leq e^{-c_5 n}$.

Then, by standard epsilon net argument in [7, p.94], we know that there exists a constant $c_6 > 0$ such that

$$\mathbb{P}(\{\omega \in \Omega : \exists x (x \in Comp(\delta, \delta n^{1/2-\gamma}), A'(\omega)x = 0, \|A'(\omega)\| \leq c_4 n^\gamma)\}) < e^{-c_6 n}$$

Therefore, we obtain that

$$\mathbb{P}(\text{LCD}_L(\frac{Z_{[1:k(n,Z)]}}{\|Z_{[1:k(n,Z)]}\|_2}) < \exp(c_{\text{Le } 6.8.1}n^{2-2\gamma}), Z \in Comp(\delta, \delta n^{1/2-\gamma}))$$
$$\leq \mathbb{P}(Z \in Comp(\delta, \delta n^{1/2-\gamma}))$$
$$\leq e^{-c_7 n}$$

for some constant $c_7 > 0$.

Combining the estimates for the case where $Z$ is compressible and the case where $Z$ is incompressible, we conclude that

$$\mathbb{P}(\text{LCD}_L(\frac{Z_{[1:k(n,Z)]}}{\|Z_{[1:k(n,Z)]}\|_2}) < \exp(c_{\text{Le } 6.8.1}n^{2-2\gamma}))$$
$$\leq \exp(-c_9 n)$$

where $c_9 > 0$ is a constant.

Let $c_{\text{Le } 6.9.1} = c_{\text{Le } 6.8.1}$ and $c_{\text{Le } 6.9.2} = c_9$. We then have the claimed result

$$\mathbb{P}(\text{LCD}_L(\frac{Z_{[1:k(n,Z)]}}{\|Z_{[1:k(n,Z)]}\|_2}) < \exp(c_{\text{Le } 6.9.1}n^{2-2\gamma})) < \exp(-c_{\text{Le } 6.9.2}n)$$

$\square$

Now, since the truncated LCD of $Z$ is large, we have a good estimate of $\mathbb{P}(|\langle Z, X_n \rangle| < \epsilon)$, and then we can use it to get the result on the invertibility for incompressible vectors.

**Theorem 6.10** (Invertibility for incompressible vectors). *Let $R_n(\omega)$ be an $n \times n$ matrix whose entries are independent copies of a centered subgaussian real random variable of unit variance. Let $M$ be an arbitrary deterministic matrix such that $\|M\| \leqslant n^\gamma$. Let $A_n(\omega) = R_n(\omega) + M$. Then there exists constants $c_{\text{Th}\,6.10.1} > 0, c_{\text{Th}\,6.10.2} > 0, c_{\text{Th}\,6.10.3} > 0, n_{\text{Th}\,6.10.1} \in \mathbb{N}$ such that for any $\gamma$ with $1/2 < \gamma < 1$, and $n > n_{\text{Th}\,6.10.1}$, we have*

$$\mathbb{P}(\inf_{x \in Incomp(\delta, \delta n^{1/2-\gamma})} \|A_n x\|_2 \leqslant \epsilon n^{-\gamma}) \leq c_{\text{Th}\,6.10.2} n^{\gamma - 1/2} (\epsilon + \exp(-c_{\text{Th}\,6.10.1} n^{2-2\gamma}))$$

*Proof.* The proof is based on the previous lemma and [8, Lemma 3.5 (Invertibility via distance), p.615].

Let $X_1(\omega), ..., X_n(\omega)$ denote the column vectors of $A_n(\omega) = R_n(\omega) + M$, and let $H_k$ denote the span of all column vectors except the $k$th. Then, [8, Lemma 3.5 (Invertibility via distance), p.615], asserts that

$$\mathbb{P}(\inf_{x \in Incomp(\delta, \rho)} \|A_n x\|_2 \leqslant \epsilon \rho n^{1/2}) \leq \frac{1}{\delta n} (\sum_{l=1}^n \mathbb{P}(\text{dist}(X_l, H_l) < \epsilon))$$

In our setting, we take the incompressibility parameter $\rho$ to be $\delta n^{1/2-\gamma}$, so we obtain

$$\mathbb{P}(\inf_{x \in Incomp(\delta, \delta n^{1/2-\gamma})} \|A_n x\|_2 \leqslant \epsilon \delta n^{-\gamma}) \leq \frac{1}{\delta n} (\sum_{l=1}^n \mathbb{P}(\text{dist}(X_l, H_l) < \epsilon))$$

So, we need to estimate $\mathbb{P}(\text{dist}(X_l, H_l) < \epsilon)$ for $l = 1, 2, ..., n$. Since we will get a universal estimate not depending on $l$, let us just work with $\mathbb{P}(\text{dist}(X_n, H_n) < \epsilon)$.

Let $Z(\omega)$ be a random unit vector such that it is orthogonal to $\text{span}(\{X_1(\omega), ..., X_{n-1}(\omega)\})$ for every $\omega \in \Omega$.

With the observation that $\text{dist}(X_n, H_n) \geq |\langle Z, X_n \rangle|$, we obtain

$$\mathbb{P}(\text{dist}(X_n, H_n) < \epsilon) \leq \mathbb{P}(|\langle Z, X_n \rangle| < \epsilon)$$

Recall that $Z(\omega)$ is a function of $X_1(\omega), ..., X_{n-1}(\omega)$. So $Z(\omega)$ is independent with $X_n(\omega)$. Let $p(z, \epsilon) = \mathbb{P}(|\langle z, X_n \rangle| < \epsilon)$, then by Fubini, we have $\mathbb{P}(|\langle Z, X_n \rangle| < \epsilon) = \mathbb{E}(p(Z(\omega), \epsilon))$.

Let $L$ be some fixed number. According to Lemma 6.9, there exists constants $c_{\text{Le}\,6.9.1} > 0, c_{\text{Le}\,6.9.2} > 0, n_1 \in \mathbb{N}$, and a function $k(n, v)$, such that for any $\gamma$ with $1/2 < \gamma < 1$, and $n > n_1$, we have $\mathbb{P}(\text{LCD}_L(\frac{Z_{[1:k(n,Z)]}}{\|Z_{[1:k(n,Z)]}\|_2}) < \exp(c_{\text{Le}\,6.9.1} n^{2-2\gamma})) < \exp(-c_{\text{Le}\,6.9.2} n)$.

We start from the following simple observation

$$\mathbb{E}(p(Z(\omega), \varepsilon))$$
$$= \mathbb{E}(1_{\{\omega \in \Omega : \text{LCD}_L(\frac{Z_{[1:k(n,Z)]}}{\|Z_{[1:k(n,Z)]}\|_2}) < \exp(c_{\text{Le}\,6.9.1} n^{2-2\gamma})\}} p(Z(\omega), \varepsilon))$$
$$+ \mathbb{E}(1_{\{\omega \in \Omega : \text{LCD}_L(\frac{Z_{[1:k(n,Z)]}}{\|Z_{[1:k(n,Z)]}\|_2}) \geq \exp(c_{\text{Le}\,6.9.1} n^{2-2\gamma})\}} p(Z(\omega), \varepsilon))$$

We first estimate the first term. By Lemma 6.9, since $p(Z(\omega), \varepsilon) \leq 1$, we have

$$\mathbb{E}(1_{\{\omega \in \Omega : \text{LCD}_L(\frac{Z_{[1:k(n,Z)]}}{\|Z_{[1:k(n,Z)]}\|_2}) < \exp(c_{\text{Le}\,6.9.1} n^{2-2\gamma})\}} p(Z(\omega), \varepsilon))$$
$$< 1 \cdot \mathbb{P}(\{\omega \in \Omega : \text{LCD}_L(\frac{Z_{[1:k(n,Z)]}}{\|Z_{[1:k(n,Z)]}\|_2}) < \exp(c_{\text{Le}\,6.9.1} n^{2-2\gamma})\})$$
$$< \exp(-c_{\text{Le}\,6.9.2} n)$$

For the second term, we use Lemma 3.5 or [10, Corollary 7.6 (Small ball probabilities for sums).,p.1742], and yields

$$p(z,\varepsilon) \leqslant c_3 n^{\gamma-1/2} L(\varepsilon + \frac{1}{\text{LCD}_L(\frac{Z_{[1:k(n,z)]}}{\|Z_{[1:k(n,z)]}\|_2})})$$

for some constant $c_3$.

Therefore, on the event $\{\omega \in \Omega : \text{LCD}_L(\frac{Z_{[1:k(n,Z)]}}{\|Z_{[1:k(n,Z)]}\|_2}) \geq \exp(c_{\text{Le }6.9.1}n^{2-2\gamma})\}$, we are able to show that

$$p(Z(\omega),\varepsilon)$$
$$\leq c_3 n^{\gamma-1/2} L(\varepsilon + \frac{1}{\text{LCD}_L(\frac{Z_{[1:k(n,Z)]}}{\|Z_{[1:k(n,Z)]}\|_2})})$$
$$\leq c_3 n^{\gamma-1/2} L(\varepsilon + \exp(-c_{\text{Le }6.9.1}n^{2-2\gamma}))$$

Therefore, we have

$$\mathbb{E}(1_{\{\omega\in\Omega:\text{LCD}_L(\frac{Z_{[1:k(n,Z)]}}{\|Z_{[1:k(n,Z)]}\|_2})\geq\exp(c_{\text{Le }6.9.1}n^{2-2\gamma})\}} p(Z(\omega),\varepsilon))$$
$$\leq \mathbb{E}(c_3 n^{\gamma-1/2} L(\varepsilon + \exp(-c_{\text{Le }6.9.1}n^{2-2\gamma})))$$
$$\leq c_3 n^{\gamma-1/2} L(\varepsilon + \exp(-c_{\text{Le }6.9.1}n^{2-2\gamma}))$$

Therefore, we derive that

$$\mathbb{E}(p(Z(\omega),\varepsilon)) \leq \exp(-c_{\text{Le }6.9.2}n) + c_3 n^{\gamma-1/2} L(\varepsilon + \exp(-c_{\text{Le }6.9.1}n^{2-2\gamma}))$$

Therefore, combining the estimates for the first term and for the second term, we have

$$\mathbb{P}(\text{dist}(X_n, H_n) < \epsilon)$$
$$\leq \mathbb{P}(|\langle Z, X_n \rangle| < \epsilon)$$
$$\leq \exp(-c_{\text{Le }6.9.2}n) + c_3 n^{\gamma-1/2} L(\varepsilon + \exp(-c_{\text{Le }6.9.1}n^{2-2\gamma}))$$
$$\leq c_4 n^{\gamma-1/2} L(\varepsilon + \exp(-c_{\text{Le }6.9.1}n^{2-2\gamma}))$$

where $c_4 > 0$ is some constant.

Since the estimate is universal, it follows that

$$\mathbb{P}(\text{dist}(X_l, H_l) < \epsilon) \leq c_4 n^{\gamma-1/2} L(\varepsilon + \exp(-c_{\text{Le }6.9.1}n^{2-2\gamma}))$$

for any $l = 1, 2, ..., n$.

Therefore, we have

$$\mathbb{P}(\inf_{x \in \text{Incomp}(\delta,\delta n^{-1/2})} \|A_n x\|_2 \leqslant \epsilon\delta n^{-\gamma}) \leq \frac{1}{c}(c_4 n^{\gamma-1/2} L(\varepsilon + \exp(-c_{\text{Le }6.9.1}n^{2-2\gamma})))$$

which completes the proof of the theorem. $\square$

## 7. Invertibility on the Whole Sphere

Combining the results for compressible and incompressible vectors, we can now complete the proof of the result of invertibility on the whole sphere, which is our main result, Theorem 1.2.

*Proof of Theorem 1.2.* The result follows directly from Lemma 5.1 and Theorem 6.10.

By Lemma 5.1, there exists $c_{\text{Le}\,5.1.1} > 0, c_{\text{Le}\,5.1.2} > 0, c_{\text{Le}\,5.1.3} > 0$ that depends only on $K$, and such that, for all $0 < \delta < c_{\text{Le}\,5.1.1}$, we have

$$\mathbb{P}(\inf_{x \in Comp(\delta, \delta n^{1/2-\gamma})} \|A_n x\|_2 \leqslant c_{\text{Le}\,5.1.2} n^{1/2}, \text{ and } \|R_n\| \leqslant K n^{1/2}) \leq e^{-c_{\text{Le}\,5.1.3} n}$$

Also, it follows from [7, p.89] that there exists a constant $c_1 > 0$ such that $\mathbb{P}(\|R_n(\omega)\| > Kn^{1/2}) < \exp(-c_1 n)$.

Hence, there exists constant $c_2 > 0$ such that

$$\mathbb{P}(\inf_{x \in Comp(\delta, \delta n^{1/2-\gamma})} \|A_n x\|_2 \leqslant c_{\text{Le}\,5.1.2} n^{1/2}) \leq \exp(-c_2 n)$$

As [8, p.614] suggests, we observe that this estimate for the class of compressible vectors is stronger than what we need because for $n$ large enough and $\epsilon$ small enough, we always have $c_{\text{Le}\,5.1.2} n^{1/2} \geq \epsilon n^{-\gamma}$.

Combining the above observation with Theorem 6.10, we obtain the desired result of invertibility on the whole sphere. $\square$

Department of Mathematics, University of Michigan, Ann Arbor, Michigan
*E-mail address*: xydong@umich.edu